\documentclass[11pt,twoside]{amsart}
\usepackage{bbm}
\usepackage{amsmath}
\usepackage{amsthm}
\usepackage{amssymb}
\usepackage{hyperref}
\usepackage[all]{xy}
\usepackage[latin1]{inputenc}

\newtheorem*{Thm}{Theorem}
\newtheorem*{Prop}{Proposition}
\newenvironment{Proof}{\par \noindent \textit{Proof.}}{\hfill $\Box$ \par}

\newcommand{\Ab}{\mathcal{A}}
\newcommand{\OO}{\mathcal{O}}
\newcommand{\EE}{\mathcal{E}}
\newcommand{\LL}{\mathcal{L}}

\newcommand{\FF}{\mathcal{F}}

\newcommand{\NNN}{\mathcal{N}}

\newcommand{\QQQ}{\mathcal{Q}}

\newcommand{\RR}{\mathbb{R}}
\newcommand{\NN}{\mathbb{N}}
\newcommand{\QQ}{\mathbb{Q}}
\newcommand{\PP}{\mathbb{P}}

\DeclareMathOperator{\rk}{rk}

\DeclareMathOperator{\Sym}{Sym}

\DeclareMathOperator{\rank}{rk}

\DeclareMathOperator{\Hilb}{Hilb}
\DeclareMathOperator{\ch}{ch}
\DeclareMathOperator{\cc}{c}

\begin{document}
\title{Stability of tautological vector bundles on Hilbert squares of surfaces}
\author{Ulrich Schlickewei}
\address{Mathematisches Institut der Universit\"at Bonn, Endenicher Allee 60, 53115 Bonn, Germany}
\email{uli@math.uni-bonn.de}

\begin{abstract}
We prove stability of rank two tautological bundles on the Hilbert square of a surface (under a mild 
positivity condition) and compute their Chern classes.
\end{abstract}

\maketitle
\let\thefootnote\relax\footnotetext{This work was supported by the SFB/TR 45 `Periods,
Moduli Spaces and Arithmetic of Algebraic Varieties' of the DFG
(German Research Foundation) and by the Bonn International Graduate School in Mathematics (BIGS).}

 Let $S$ be a smooth, projective surface, let $\Hilb^2(S)$ be the Hilbert scheme
 parametrizing subschemes of $S$ of length 2. It is known by a classical theorem
 of Fogarty \cite{Fo} that $\Hilb^2(S)$ is a smooth, projective variety of dimension 4.
 Let $Z \subset S \times \Hilb^2(S)$ be the universal subscheme, denote by $p: Z \to S$ and
 by $q: Z \to \Hilb^2(S)$ the projections. Given a line bundle $\LL$ on $S$, the sheaf
 $\LL^{[2]} := q_* p^* \LL$ is a rank two vector bundle on $\Hilb^2(S)$, called the \emph{tautological vector bundle associated with
 $\LL$}. 
 
 In this note we prove the following
 
    \begin{Thm} Assume that $h^0(S,\LL)\ge 2$. 
     Then for $N \gg 0$, the vector bundle $\LL^{[2]}$ is $\mu_{H_N}$-stable
     on $\Hilb^2(S)$.
     \end{Thm}
     
     Here, $H_N$ is a polarization of the form $\Sym^2(NH) - E$ where $H$ 
     is an ample divisor on $S$ and $E \subset \Hilb^2(S)$ denotes the
     exceptional divisor of the Hilbert--Chow morphism.
     
     The proof of the theorem relies upon the fundamental short exact sequence 
     for tautological vector bundles on the blowup of $S \times S$ and
     upon the corresponding result for curves which was proved by Mistretta \cite{Mi}.

     \vspace{1ex}
     Originally, our interest in this result came from the desire to produce 
     vector bundles on Hilbert schemes of
     K3 surfaces with interesting metrics and with interesting Chern classes.
     For this reason we give a formula for the Chern
     classes of $\LL^{[2]}$ in terms of the symmetric product of $\cc_1(\LL)$, of $[E]$ and of the
     characteristic classes of $\Hilb^2(S)$.
     
     \vspace{1ex}     
     After introducing some notation in Section \ref{notation} we proof the theorem in Section \ref{Proof}.
     Finally we calculate the Chern classes of $\LL^{[2]}$ in Section \ref{Chern}.

    \vspace{1ex}
    \noindent
   {\emph{Acknowledgements.} This paper is a part of my Ph.D.\ thesis prepared
    at the University of Bonn. It is a great pleasure to thank my advisor Daniel Huybrechts for 
    his constant support.
    I am also grateful to
    Luca Scala for helpful discussions on Chern characters
    of tautological vector bundles and to 
		Ernesto Mistretta for explaining to me his results about stable vector bundles on
   	symmetric products of curves.

     \begin{section}{Some notation} \label{notation}
     Let $\iota_{\Delta}: \Delta \hookrightarrow S \times S$ be the diagonal.
     Denote by $\sigma: \widetilde{S \times S} \to S \times S$ the blowup of $S \times S$ in 
$\Delta$. The natural action of the symmetric group $\mathfrak{S}_2$ on $S \times S$ extends
to a holomorphic action on $\widetilde{S \times S}$ and $\Hilb^2(S) = \widetilde{S \times S} / 
\mathfrak{S}_2$. Let $\iota_D : D \hookrightarrow \widetilde{S \times S}$ be the 
exceptional divisor of $\sigma$.

In the following diagram we summarize the situation and, at the same time,
give names to the various natural maps.
 \begin{equation*}
  \xymatrix{
   D \ar[rr]^{\iota_D} \ar[d]_{\sigma_D} & & \widetilde{S \times S} \ar@/_/[ddl]_{\widetilde{r}_1} 
   \ar@/^/[ddr]^{\widetilde{r}_2} \ar[d]_{\sigma} \ar[rr]^-{\pi} & & \Hilb^2(S) 
   \ar[d] \\
   S \simeq \Delta \ar[rr]_{\iota_{\Delta}} & &  S \times S \ar[dl]^{r_1} \ar[dr]_{r_2} \ar[rr] & & \Sym^2(S) \\
   & S & & S}
  \end{equation*}
  
  Given an ample hypersurface $H \subset S$, for $N$ sufficiently large the divisor $\widetilde{H}_N :=
  N(\widetilde{r}_1^*H + \widetilde{r}_2^*H) - D$ is ample on $\widetilde{S \times S}$. Moreover, this divisor
  is invariant under the action of $\mathfrak{S}_2$ on the divisor class group of $\widetilde{S \times S}$, hence 
  it is of the form $\pi^* H_N$ for an ample divisor $H_N$ on $\Hilb^2(S)$.
  
  \vspace{1ex}
  
  Finally we recall the notion of $\mu$-stability. Let $Y$ be a smooth, projective variety, polarized by an
  ample divisor $H$. Let $\EE$ be a torsion-free coherent $\OO_Y$-module. Then the \emph{slope of $\EE$
  with respect to $H$} is defined as
  \begin{equation} \label{slope}
    \mu_H (\EE) := \frac{\int_Y \cc_1(\EE) [H]^{\dim(Y)-1}}{\rk(\EE)}.
  \end{equation}
  The sheaf $\EE$ is called \emph{$\mu_H$-stable} if for any subsheaf $\FF \subset \EE$ with 
  $0 < \rk(\FF) < \rk(\EE)$ we have $\mu_H(\FF) < \mu_H(\EE)$.
  \end{section}
 	   
 	   \begin{section}{Proof of the Theorem} \label{Proof}
      If $\LL^{[2]}$ had a destabilizing subsheaf, then by passing to the reflexive hull we see that there 
      would exist a destabilizing
      line bundle. Thus, any
      destabilizing subsheaf of $\LL^{[2]}$ on $\Hilb^2(S)$ with respect to $H_N$  induces 
      a destabilizing sub-line bundle of $\EE:= \pi^* \LL^{[2]}$ on $\widetilde{S \times S}$ 
      with respect to $\widetilde{H}_N$. We will show that $\EE$ is
      $\mu_{\widetilde{H}_N}$-stable.
      This will finish the proof of the theorem.
      
      \vspace{1ex}
      For $i=1,2$, put $\LL_i := \widetilde{r}_i^* \LL$ and let $\LL_D := \iota_{D,*} \sigma_D^* \LL$. 
      Consider the fundamental short exact 
      sequence (see e.g.\ the proof of \cite[Prop.\ 2.3]{Da})
      \begin{equation} \label{SES_E}
        0 \to \EE \to \LL_1 \oplus \LL_2 \to \LL_D \to 0
      \end{equation}
      where the surjection on the right is given by $(s_1, s_2) \mapsto s_{1|D} - s_{2|D}$. 
      This sequence shows that $\cc_1(\EE) = \cc_1(\LL_1) + \cc_1(\LL_2) - [D]$.
      Moreover, we deduce that 
      $\EE$ contains the line bundles $\LL_1(-D)$ and $\LL_2(-D)$.

      \vspace{1ex}
      Let now $\Ab \subset \EE$ be an arbitrary sub-line bundle. Then $\Ab$ has one of the following three
      properties:
      
      1.) $\Ab \subset \LL_1(-D)$,
      
      2.) $\Ab \subset \LL_2(-D)$,
      
      3.) $\Ab \not\subset \LL_1(-D)$ and $\Ab \not\subset \LL_2(-D)$.
      
      We will prove that there exist $N_1, N_2, N_3 \in \NN$ such that for all $\Ab \subset \EE$ with
      property $i = 1,2$ or 3 and for all $N \ge N_i$ we have
      \begin{equation} \label{inequality_N_i}
        \mu_{\widetilde{H}_N} (\Ab) < \mu_{\widetilde{H}_N} (\EE).
      \end{equation}

      \vspace{1ex}
      Assume first that we are in case 1.), i.e.\ that $\Ab \subset \LL_1(-D)$. Choose
      a natural number $N_1 \ge 4$ (this will be useful in (\ref{N1>4}) below) such that
      $\widetilde{H}_N$ is ample for all $N \ge N_1$. Then 
      $\mu_{\widetilde{H}_N} (\Ab) \le \mu_{\widetilde{H}_N}(\LL_1(-D))$ for all $N \ge N_1$.
      Let $\alpha_i := \cc_1(\LL_i)$. Then $\cc_1 (\EE) = 
      \alpha_1 + \alpha_2 - [D]$ and therefore
      
      \begin{equation} \begin{aligned}
        \mu_{\widetilde{H}_N} (\EE) - \mu_{\widetilde{H}_N}(\Ab) \ge & \; 
        \mu_{\widetilde{H}_N} (\EE) - \mu_{\widetilde{H}_N}(\LL_1(-D)) \\
        = & \; \int_{\widetilde{S \times S}} \left( \frac{\alpha_1 + \alpha_2 - [D]}{2} - (\alpha_1 - [D]) \right) 
         [\widetilde{H}_N]^3 \\
        = & \; \int_{\widetilde{S \times S}} 
              \frac{\alpha_2 - \alpha_1}{2} [\widetilde{H}_N]^3 + 
              \int_{\widetilde{S \times S}} \frac{[D]}{2} [\widetilde{H}_N]^3 \\
        = & \; 0 + \int_{\widetilde{S \times S}} \frac{[D]}{2} [\widetilde{H}_N]^3 \\
        > & \; 0
        \end{aligned}
      \end{equation}
      for all $N \ge N_1$ because $D$ is effective and $\widetilde{H}_N$ is ample.
       
      \vspace{1ex}
      An analogous reasoning applies in case 2.) with $N_2 = N_1$.
      
      \vspace{1ex}
      In case 3.) we proceed in two steps. First we show that
      \begin{equation} \label{inequality_mu_F}
        \mu_F(\Ab) < \mu_F (\EE),
      \end{equation}
      where $F = \widetilde{r}_1^* H + \widetilde{r}_2^* H$ and $\mu_F$ is the slope with
      respect to the nef divisor $F$, defined as in (\ref{slope}). 
      Then we use an asymptotic argument to complete
      the proof.
      
      \vspace{1ex}
      To prove (\ref{inequality_mu_F}), we 
      will consider two divisors in $|MF|$ for an appropriate $M>0$
      which intersect along a reducible surface.
      Then we reduce our computation to the irreducible components of these surfaces.

      Choose $M$ sufficiently positive such that the linear system $|MH|$ 
      contains two distinct, smooth curves $C$ and $C'$ meeting
      transversely. Denote by
      \begin{equation*}
        G_1:= (C \times S) \cup (S \times C') \; \; \text{and by} \; \; G_2 := (C' \times S) \cup (S \times C),
      \end{equation*}
      let $\widetilde{G}_1$ and $\widetilde{G}_2$ be their strict transforms under $\sigma: \widetilde{S \times S} \to S \times S$.
      Note that $\widetilde{G}_i \in | MF |$ because each component
      of $G_i$ meets the center of the blowup, namely the diagonal of $S \times S$, along a curve. 
      The intersection $\widetilde{G}_1 \cap \widetilde{G}_2$ is the disjoint
      union of the four smooth surfaces
      \begin{equation*} \begin{aligned}
        T_1 : = \widetilde{(C \cap C') \times S}, \; \; & \; \; T_2:= \widetilde{S \times (C \cap C')}, \\
        T_3:= \widetilde{C \times C}, \; \; & \; \; T_4:= \widetilde{ C' \times C' }, \end{aligned}
      \end{equation*}
      where for any subvariety $Y \subsetneq S \times S$ we write $\widetilde{Y}$ for the strict transform of $Y$ under
      $\sigma: \widetilde{S \times S} \to S \times S$.
      
      Then for the $F$-slope of any coherent sheaf $\mathcal{F}$ we find
      \begin{equation*}
        \mu_{MF} (\mathcal{F}) = \frac{\deg_{MF} (\mathcal{F})}{\mathrm{rk}(\mathcal{F})} = \sum_{i=1}^4 
        \frac{\deg_{T_i} (\mathcal{F})}{\mathrm{rk}(\mathcal{F})}
      \end{equation*}
      where $\deg_{T_i} (\mathcal{F}) := \int_{T_i} \cc_1(\mathcal{F}_{|T_i}) [MF]_{|T_i}$. 
      We will show that 
      \begin{equation*}
        \deg_{T_i} (\Ab) \le \frac{\deg_{T_i} (\EE)}{2}
      \end{equation*}
      for $i=1, \ldots ,4$ with strict inequality for $i= 3,4$. This will conclude the
      proof of (\ref{inequality_mu_F}).
 
      \vspace{1ex} $i=1$: The surface $T_1$ is a disjoint union of surfaces of the form $S_p := \widetilde{ \{ p \} \times S}$, 
      $p\in S$ running over the finite set $C \cap C'$.
      Since
      the fundamental class of $S_p$ does not vary for different $p \in S$, we fix an arbitrary 
      point $p_0 \in S$ and we 
      get for any coherent sheaf $\FF$
      on $\widetilde{S \times S}$
      \begin{equation*}
       \deg_{T_1}(\FF) = \sharp(C \cap C') \deg_{S_{p_0}} (\FF).
      \end{equation*}
      Note that $S_{p_0}$ is isomorphic to the blow-up of $S$ in $p_0$.
      Denote by $\sigma_{p_0}: S_{p_0} \to S$ the blow-down and by 
      $E_{p_0} \subset S_{p_0}$ the exceptional divisor. Then
      \begin{equation*}
        \cc_1 (\EE_{|S_{p_0}}) = \sigma_{p_0}^* \cc_1(\LL) - [E_{p_0}]
      \end{equation*}
      because $\cc_1(\EE) = \widetilde{r}_1^* \cc_1(\LL) + \widetilde{r}_2^* \cc_1(\LL) - [D]$ and 
      because $\widetilde{r}_1^* \LL_{| S_{p_0}} = \OO_{S_{p_0}}$.
      Now suppose that 
      \begin{equation} \label{DestabilizingT_1}
        2 \deg_{T_1}(\Ab) > \deg_{T_1} (\pi^* \LL^{[2]}).
      \end{equation}
      Then we would get
      \begin{equation*} 
        2 \int_{S_{p_0}} \cc_1 (\Ab_{|S_{p_0}})  [MF]_{|S_{p_0}} >  \int_{S_{p_0}} \big( \sigma_{p_0}^* \cc_1(\LL) - [E_{p_0}] \big) [MF]_{|S_{p_0}} 
      \ge  0
      \end{equation*}
      because $[MF]_{|S_{p_0}}$ is a nef
      class on $S_{p_0}$ and
      because $\sigma_{p_0}^* \LL \otimes \OO(-E_{p_0})$ is  the line bundle of an 
      effective divisor on $S_{p_0}$. Indeed, since $h^0(\LL) \ge 2$, there exists
      a divisor $K \in |\LL|$ with ${p_0} \in \mathrm{supp}(K)$.
      Then the strict transform $\widetilde{K}$
      of $K$ is in the
      linear system $|\sigma_{p_0}^* \LL \otimes \OO(- k E_{p_0})| $
      for some $k \ge 1$.
      Thus, $\widetilde{K} + (k-1)E_{p_0}$ is an effective divisor with
      line bundle 
      $\sigma_{p_0}^* \LL \otimes \OO(-E_{p_0})$. 
      
      By (\ref{SES_E}), $\Ab_{|S_{p_0}} \subset \big( \LL_1 \oplus \LL_2 \big)_{|S_{p_0}} = \OO_{S_{p_0}} \oplus \sigma_{p_0}^* \LL$.
      Now, since $\Ab_{|S_{p_0}}$ has positive $F_{|S_{p_0}}$-slope, 
      the composition $\Ab_{|S_{p_0}} \to \OO_{S_{p_0}} \oplus \sigma_{p_0}^* \LL \to \OO_{S_{p_0}}$ must
      be zero. Since we chose $p_0 \in S$ randomly, 
      assumption (\ref{DestabilizingT_1}) implies that this composition is zero for all $p \in S$.
      On the other hand, all $x \in \widetilde{S \times S}$ 
      lie on some $S_p$. This shows that the composition
      \begin{equation*}
        \Ab \to \LL_1 \oplus \LL_2 \to \LL_1
      \end{equation*}
      is zero. But then the short exact sequence (\ref{SES_E}) implies that $\Ab \subset 
      \widetilde{r}_2^* \LL(-D)$, because
      the surjection $\LL_1 \oplus \LL_2 \to \LL_D$ is given
      by $(s_1 , s_2) \mapsto s_{1|D} - s_{2|D}$. 
      This is a contradiction to the assumption that $\Ab$ satisfies 3.).

      \vspace{1ex} $ i=2$: analogous to $i=1$.

      \vspace{1ex} $ i = 3$: Note that $\widetilde{C \times C}$ is isomorphic to $C \times C$ and that 
      \begin{equation*}
        MF_{| \widetilde{C \times C}} = M (p_1^* H_{|C} + p_2^* H_{|C})
      \end{equation*}
      where $p_i: C \times C \to C$ are the projections. Moreover, it is easily checked that 
      \begin{equation*}
        (\pi^* \LL^{[2]})_{| \widetilde{C \times C}} = \pi_C^* \LL_{|C}^{[2]}
      \end{equation*}
      where $\pi_C: C \times C \to \Hilb^2(C)$ is
      the natural projection and $\LL_{|C}^{[2]}$ is the tautological line bundle associated with $\LL_{|C}$ on $\Hilb^2(C)$.
      It remains to apply \cite{Mi}, Cor.\ 4.3.3 which says that $\pi_C^* \LL_{|C}^{[2]}$ is a stable vector bundle on $C \times C$.

      \vspace{1ex} $i=4$: analogous to $i=3$. 
      
      \vspace{1ex} Thus, we have proved (\ref{inequality_mu_F}). 
      To conclude the proof of (\ref{inequality_N_i}) we define for $n \in \NN$ 
      the linear function
      \begin{equation*} 
        \varphi_n: K^0(\mathrm{Coh}(\widetilde{S \times S})) \to \QQ, \; \; \;
        \mathcal{F} \mapsto \sum_{i=0}^2 n^i \frac{3 !} {i! (3-i)!}
        \int_{\widetilde{S \times S}}  \cc_1(\mathcal{F}) [F]^i [\widetilde{H}_{N_1}]^{3-i}.
      \end{equation*}
      Then noting that $\widetilde{H}_{n + N_1} = n F + \widetilde{H}_{N_1}$, we get 
      for all $\mathcal{F} \in \mathrm{Coh}(\widetilde{S \times S})$
      \begin{equation} \label{mu_H=mu_F+phi_n}
        \mu_{\widetilde{H}_{n + N_1}} (\mathcal{F}) = n^3 \mu_{F} (\mathcal{F}) + 
        \frac{\varphi_n (\mathcal{F})}{\rank(\mathcal{F})}.
      \end{equation} 
      Inequality (\ref{inequality_mu_F}) implies that there exists a positive constant $k \in \RR_{> 0}$ such that
      for all sub-line bundles $\Ab \subset \EE$ with property 3.) 
      we have 
      \begin{equation} \label{mu_F_k}
        \mu_F (\EE) - \mu_F(\Ab) \ge k.
      \end{equation}
      This is because 
      $\mu_F$ takes integer values on line bundles.
      
      We will now show that $\varphi_n(\Ab) < \varphi_n (\EE)$. To see this, we first prove that
      $\Ab \cap \LL_1(-D) = 0$. Indeed, otherwise the torsion-free sheaf $\Ab + \LL_1(-D)$ would be of 
      rank 1 as follows from the short exact sequence
      \begin{equation*}
       0 \to \Ab \cap \LL_1(-D) \to \Ab \oplus \LL_1(-D) \to \Ab + \LL_1(-D) \to 0.
      \end{equation*}
      Then the reflexive hull $\Ab'$ of $\Ab + \LL_1(-D)$ would be a sub-line bundle of $\EE$ which 
      again would have property 3.) because $\Ab \subset \Ab'$. From (\ref{inequality_mu_F}) we deduce
      $\mu_F(\Ab') < \mu_F(\EE)$. 
      On the other hand 
      \begin{equation*} \begin{aligned}
        \mu_F (\LL_1(-D)) = & \; \int_{\widetilde{S \times S}} (\alpha_1 - [D]) [F]^3 \\
          = & \; \int_{\widetilde{S \times S}} \alpha_1 [F]^3 \\
          = & \; \int_{\widetilde{S \times S}} \frac{\alpha_1 + \alpha_2}{2} [F]^3 \\
          = & \; \int_{\widetilde{S \times S}} \frac{\alpha_1 + \alpha_2 - [D]}{2} [F]^3 \\
          = & \; \mu_F (\EE).
          \end{aligned}
      \end{equation*}
      Here, we used that the integral over a cohomology class of degree $6$ on $D$ which
      is pulled back from $\Delta$ vanishes.
      Using that $\LL_1(-D) \subset \Ab'$ and that $F$ is nef, we get a chain of inequalities
      \begin{equation*}
        \mu_F(\EE) = \mu_F(\LL_1(-D)) \le \mu_F (\Ab') < \mu_F(\EE).
      \end{equation*}   
      This is a contradiction, whence $\Ab \cap \LL_1(-D) = 0$. 
      
      Then we have a short exact sequence
      \begin{equation} \label{SES_Q}
        0 \to \Ab \oplus \LL_1(-D) \to \EE \to \QQQ \to 0
      \end{equation}
      where $\QQQ$ is a torsion sheaf. It follows that $\cc_1(\QQQ)$ is either zero or effective.
      Since $\varphi_n$ involves only products of the nef divisor $F$ and the ample divisor $\widetilde{H}_{N_1}$,
      this implies that $\varphi_n(\QQQ) \ge 0$. 
      
      We claim that
      there exists $n_1 \in \NN$ such that
      for $n \ge n_1$ we have $\varphi_n(\LL_1(-D)) >0$.
      To see this, it is enough to show that the $n^2$-term of $\varphi_n(\LL_1(-D))$ is
      positive. We have
      \begin{equation*} \begin{aligned}
        \cc_1(\LL_1(-D)) [F]^2 [\widetilde{H}_{N_1}]  = & \; (\alpha_1 - [D]) [F]^2 (N_1 [F] -[D]) \\
        = & \; N_1 \alpha_1 [F]^3 + [F]^2 [D]^2 - [D]( N_1 [F]^3 + \alpha_1 [F]^2).
        \end{aligned}
      \end{equation*}
      If $q$ denotes
      the intersection product on $S$, then we obtain
      \begin{equation} \label{N1>4} \begin{aligned}
        \int_{\widetilde{S \times S}} N_1 \alpha_1 [F]^3 + \int_{\widetilde{S \times S}} [F]^2 [D]^2 = & \; 
        N_1 \int_{S \times S} r_1^* \cc_1 (\LL) \left( r_1^* [H] + r_2^* [H] \right)^3 + \int_D [F]^2 \xi \\
        = & \; \bigg( 3 N_1 q \big(\cc_1 (\LL), [H] \big) - 4  \bigg) q([H]) \\
        > & \; 0. \end{aligned}
      \end{equation}
      Here, $r_i : S \times S \to S$ are the projections and $\xi = [D]_{|D} = 
      \cc_1 (\OO_{\PP(\NNN_{\Delta | S \times S})}(-1))$.
      We use that $\int_D (\sigma_D^* \alpha) \xi = - \int_S \alpha$ for $\alpha \in H^*(S,\QQ)$,
      that $N_1 \ge 4$ and that $\cc_1(\LL)$ is an effective
      class on $S$. This proves the existence of $n_1$ with $\varphi_n (\LL_1(-D)) > 0$
      for all $n \ge n_1$.  
         
      Now by (\ref{SES_Q}) we have
      \begin{equation} \label{phi_n(A)}
        \varphi_n(\Ab) = \varphi_n (\EE) - \varphi_n(\QQQ) - \varphi_n (\LL_1(-D)) < \varphi_n(\EE)
      \end{equation}
      for all $n \ge n_1$.
      
      Putting together (\ref{mu_H=mu_F+phi_n}), (\ref{mu_F_k}) and (\ref{phi_n(A)}) we find for
      $n \ge n_1$ and for all line bundles $\Ab \subset \EE$ with property 3.)
      \begin{equation*} \begin{aligned}
        \mu_{\widetilde{H}_{n + N_1}}(\EE) - \mu_{\widetilde{H}_{n + N_1}} (\Ab) = & \; 
        n^3 (\mu_F(\EE) - \mu_F(\Ab)) + \frac{\varphi_n(\EE)}{2} - \varphi_n(\Ab) \\
        > & \; n^3 k - \frac{\varphi_n (\EE)}{2}.
      \end{aligned} \end{equation*}
      Now since $k >0$ and since $\varphi_n(\EE)$ is a polynomial of degree 2 in $n$,
      there exists $n_2 \ge n_1$ such that $n^3 k - \frac{\varphi_n(\EE)}{2} >0$ for
      all $n \ge n_2$. Therefore,
      with $N_3 := N_1 + n_2$, inequality (\ref{inequality_N_i}) is satisfied for $i=3$. This completes
      the proof.      \qed
      \end{section} 
      
      \begin{section}{The Chern character} \label{Chern}
      In this section we express the Chern classes of $\LL^{[2]}$ in terms of $\cc_2(\Hilb^2(S))$, of
      the symmetric product of $\cc_1(\LL)$ and of the fundamental class the the exceptional divisor
      of the Hilbert--Chow morphism $\Hilb^2(S) \to \Sym^2(S)$. 
      Since $\Hilb^2(S) \simeq \widetilde{S \times S} / \mathfrak{S}_2$, via pullback along the quotient
      morphism $\pi$ we
      get an identification $H^*(\Hilb^2(S),\QQ) \simeq H^*(\widetilde{S \times S}, \QQ)^{\mathfrak{S}_2}$
      (see \cite{Gr}) and we are reduced to calculate the Chern classes of $\pi^* \LL^{[2]}$.
      By the short exact sequence (\ref{SES_E}) we get
 \begin{equation*} \label{Formel1Chern} \begin{aligned}
   \pi^* \cc_1(\LL^{[2]}) & = \alpha_1 + \alpha_2 - [D] \; \; \text{and} \\
   \pi^* \cc_2(\LL^{[2]}) & = \cc_2(\LL_1 \oplus
     \LL_2) - \cc_2(\LL_D) - \cc_1(\LL^{[2]}) \cc_1(\LL_D) = \alpha_1 \alpha_2 - \frac{1}{2} (\alpha_1+\alpha_2) [D],
   \end{aligned}
 \end{equation*}
 where as above $\alpha_i = \widetilde{r}_i^* \cc_1(\LL)$ and where we used the Grothendieck--Riemann--Roch theorem
 to calculate $\cc_2(\LL_D)$.
 
 Using the tangent bundle sequence for the ramified covering $\pi: \widetilde{S \times S} \to \Hilb^2(S)$
 and the formula for Chern classes of blow-ups (cf.\ \cite[Ex.\ 15.4.3]{F}) 
 we can express the pullback of the $(4,0)+(0,4)$-K\"unneth factor of the diagonal $\Delta \subset S \times S$ 
 to $\widetilde{S \times S}$ as follows:
 \begin{equation*}
   -(\int_S \ch_2(S)) \cdot \sigma^* ([\Delta]^{4,0} + [\Delta]^{0,4}) = 
   \pi^* \big(\cc_2(\Hilb^2(S)) - \frac{\cc_1^2(\Hilb^2(S)) + \cc_1(\Hilb^2(S)) \delta}{2} + 3 \delta^2 \big)
 \end{equation*}
 where $\delta \in H^2(\Hilb^2(S),\QQ)$ is the class with $\pi^* \delta = [D]$.
 If $\int_S \ch_2(S) \neq 0$, this implies that  
 \begin{equation*} \begin{aligned}
   \alpha_1 \alpha_2 & = \frac{1}{2}\big( (\alpha_1 + \alpha_2)^2 - \alpha_1^2 - \alpha_2^2 \big) = 
      \frac{1}{2} \big( (\alpha_1 + \alpha_2)^2 - q(\cc_1(\LL)) \cdot \sigma^*([\Delta]^{4,0} + [\Delta]^{0,4}) \big) \\
   & = \frac{(\alpha_1 + \alpha_2)^2}{2} + \frac{q(\cc_1(\LL))}{2 \int_S \ch_2(S)} \big\{
       \cc_2(\Hilb^2(S)) - \frac{\cc_1^2(\Hilb^2(S)) + \cc_1(\Hilb^2(S)) \delta}{2} + 3 \delta^2 \big\}.
  \end{aligned}
 \end{equation*}
 Here, $q$ denotes the intersection product on $H^2(S,\QQ)$.

 Let $\varphi: H^2(S,\QQ) \to H^2(\Hilb^2(S),\QQ)$ be the homomorphism which is determined by
 \begin{equation*}
  \pi^* \varphi(\beta) = \widetilde{r}_1^* \beta + \widetilde{r}_2^* \beta \; \; \text{for all} \; \; \beta \in 
  H^2(S,\QQ).
 \end{equation*}
  
  Summarizing the above discussion we have achieved an expression of the Chern classes of $\LL^{[2]}$ in terms 
  of $\cc_1(\Hilb^2(S))$ and $\cc_2(\Hilb^2(S))$, of $\varphi(\cc_1(\LL))$ and of $\delta$:
  
  \begin{Prop}
    i) Assume that $\int_S \ch_2(S) \neq 0$. Then the Chern classes of $\LL^{[2]}$ are
    \begin{equation*} \begin{aligned}
      \cc_1(\LL^{[2]}) & = \varphi (\cc_1(\LL)) - \delta \; \; \text{and} \\
      \cc_2(\LL^{[2]}) & = \frac{ \varphi (\cc_1(\LL))^2 - \varphi (\cc_1(\LL)) \delta}{2} \\  
         & \; \; \; \; \; \; + \frac{q(\cc_1(\LL))}{2 \int_S \ch_2(S)} \big\{\cc_2(\Hilb^2(S)) - \frac{\cc_1^2(\Hilb^2(S)) + \cc_1(\Hilb^2(S)) \delta}{2} + 3 \delta^2 \big\}.
     \end{aligned}
    \end{equation*}
         
     ii)For $N \gg 0$, the bundle $\LL^{[2]}$ satisfies the strict Bogomolov--L\"ubke inequality, that is
     \begin{equation*}
        \Delta (\LL^{[2]}) H_N^2  = \big( 4 \cc_2(\LL^{[2]}) - \cc_1^2(\LL^{[2]}) \big) H_N^2 > 0.
     \end{equation*}
     In particular, $\LL^{[2]}$ is not projectively flat.
 \end{Prop}
 
 \begin{Proof}
  It remains to show ii). Using the above calculations we find
  \begin{equation*}
    \pi^* \Delta(\LL^{[2]}) = 2 \alpha_1 \alpha_2 - \alpha_1^2 - \alpha_2^2 - [D]^2 = 
    \Delta(\LL_1 \oplus \LL_2) - [D]^2.
  \end{equation*}
  Now, $\LL_1 \oplus \LL_2$ is $\mu_{\widetilde{H}_N}$-polystable, thus
  \begin{equation*}
    \Delta(\LL_1 \oplus \LL_2) \widetilde{H}_N^2 \ge 0.
  \end{equation*}
  On the other hand,
  \begin{equation*}
    -[D]^2 \widetilde{H}_N^2 = - \int_D \widetilde{H}_N^2 \xi = 4 N^2 q(H) + O(N) > 0 \; \; \text{for} \; N \gg 0.
  \end{equation*}
  Thus 
  \begin{equation*}
    \Delta(\LL^{[2]}) H_N^2 = \pi^* \Delta(\LL^{[2]}) \widetilde{H}_N^2 > 0
  \end{equation*}
  for $N \gg 0$.
  \end{Proof}     
   \end{section}


\begin{thebibliography}{000}

\bibitem[D]{Da} G.\ Danila, \emph{Sur la cohomologie d'un fibr\'e tautologique sur le sch\'ema de 
Hilbert d'une surface}, J.\ Alg.\ Geom.\ {\bf 10} (2001), 247-280.

\bibitem[Fo]{Fo} J.\ Fogarty, \emph{Algebraic families on an algebraic surface}, 
Am.\ J.\ Math.\ {\bf 90} (1968), 511-521.

\bibitem[Fu]{F} W.\ Fulton, \emph{Intersection theory}, Erg. Math., 3. Folge, Band 2, 2nd ed., Springer (1998).

\bibitem[G]{Gr} A.\ Grothendieck, \emph{Sur quelques points d'alg\`ebre homologique}, T\^{o}hoku Math.\ J.\ 
{\bf 9} (1957), 119-221.

\bibitem[M]{Mi} E.\ Mistretta, \emph{Some constructions around stability of vector bundles on projective varieties}, PhD Thesis, Univ.\ Paris VII (2006).
 
\end{thebibliography}
\end{document}